\input amstex
\documentstyle{amsppt}

\magnification=1200

\def\R{${\text{\bf R}}^N$}

\def\bs{\bigskip}

\def\a{\alpha}
\def\n{\nabla}

% "Sparar" check-accent som /hacek, innan /v omdefinieras.

\nologo
\NoRunningHeads

\topmatter
\title
A note on Sobolev type inequalities with $p<1$ and modified norm.
\endtitle
\author 
Andreas Wannebo
\endauthor

\endtopmatter

\document
\heading
0. Introduction
\endheading

The paper is based on one main inequality. 
Here formulated for polynomials or power series.
There ought to be other applications of the proof technique as well.
The inequality is between norm-like entities. 
The (quasi)norms used here are modifidified $L^p$-norms with
a cube as domain. The issue which domains to study is however not further
developed here.
In the inequality the upper part involves exponent $p<1$ and
its lower part the Sobolev exponent $p^*$ with
$p^*={N\cdot p\over N-1\cdot p}$ and $p^*\ge 1$.
The inequality has the form of a (weak) Poincar\'e inequality.
\par

-- The proof is easy but maybe unexpected.
\bs

\heading
1. Result
\endheading

An idea is that it is possible to treat
bad situations like this inequality for $p<1$ by denying cancellation
and this is done in the definition of the (quasi) norm.
\bs

The study is based on a ``fundamental domain'' $Q_0=(0,1)^N$
in \R. It is chosen to make every monomial have fixed
sign inside.
\bs

{\bf 1.1 Definitions.}
{\sl
Let $\Cal P$ the set of real polynomials in \R\ and
let $P_k\in \Cal P$ be the polynomials of degree less than or equal to $k$.
}
\bs

Define a ``triple stroke'' $L^p$ (quasi) norm as follows.
\bs

{\bf 1.2 Definition.} {\sl Let}

$$
|||P|||_{L^p(Q_0)}=||P_0||_{L^p(Q_0)}+||P_{\oplus}||_{L^p(Q_0)}
+||P_{\ominus}||_{L^p(Q_0)}.
$$
\bs

Here $P_0$ is the constant term of $P$.
Furthermore $P_{\oplus}$ is the sum of all its monomials with
their coefficients under condition that the degree is positive and the coefficients
are positive.
The polynomial $P_{\ominus}$ is defined in the same way
using monomials with negative coefficients instead.
\vfill\eject

Next the main result is stated. It is a weak Poincar\'e inequality
with Sobolev exponent.
\bs

{\bf 1.3 Theorem.} 
{\sl
Let $P\in\Cal P$. 
Let ${N\over N+1}\le p<1$ with $p^*={N\over N-p}$ i.e. also
$1\le p^*<{N\over N-1}$.
Let $Q_0=(0,1)^N$. 
It holds}

$$
|||P-P_0||||_{L^{p*}(Q_0)}\le A(N,p)\cdot |||\n P||||_{L^p(Q_0)}.
$$
\bs

{\bf Proof:}
Study a monomial $x^{\a}$ (multiindex notation).
First its $L^{p^*}(Q_0)$-norm is calculated,
observe that this simplifies if some $\a_j=0$,

$$
||x^{\a}||_{L^{p^*}(Q_0)}
=
\prod_1^N (\int_0^1x_i^{p^*\a_j} dx_j)^{1\over p^*}
=
\prod_1^N{1\over (p^*\a_j+1)^{1\over p*}}.
$$
\bs

Then the $p$-(quasi)norm of the 
gradient of the monomial, $\n x^{\a}$, is calculated

$$
||\n x^{\a}||_{L^p(Q_0)}=(\sum ||D_j x^{\a}||_{L^p(Q_0)}^p)^{1\over p}
=
(\sum_j \a_j^p \cdot {\a_jp+1\over \a_j-1})^{1\over p}
				\cdot
(\prod_{j=1}^N{1\over \a_jp+1})^{1\over p}.
$$
\bs

In order to establish the theorem for monomials estimate the following

$$
{||x^{\a}||_{L^{p^*}(Q_0)}\over||\n x^{\a}||_{L^p(Q_0)}}
=
\prod_1^N{(\a_jp+1)^{1\over p}
\over (p^*\a_j+1)^{1\over p^*}}
\cdot
(\sum_j \a_j^p \cdot {\a_jp+1\over \a_j-1})^{-{1\over p}}
$$

The product is estimated one factor at a time

$$
\prod_1^N{(\a_jp+1)^{1\over p}
\over (p^*\a_j+1)^{1\over p^*}}
\lesssim
\prod_1^N
\a_j^{1\over N}
\lesssim
||\a||_{l^1}.
$$

Here the coefficients are treated as a vector in a 
discrete space. The last instance is the inequality between 
geometric and arithmetic means. Furthermore the second factor is 
$
\lesssim
||\a||_{l^p}^{-1}.
$
The standard inequality between $l^q$ (quasi) norms, i.e.
they decrease when $q$ increase, implies that
the total product is less than or equal to (or similar to) 1
since $p<1$.
\bs

Hence the main inequality holds for monomials and all $p$-values
in the range.
\vfill\eject

Next step is to prove the main inequality for real polynomials
$P=\sum c_{\a} x^{\a}$.
\bs

It follows from the triangle inequality, with $1\le p^*$, that

$$
|||P|||_{L^{p^*}(Q_0)}
\le
\sum |a_{\a}|\cdot ||x^{\a}||_{L^{p^*}(Q_0)}.
$$

It can be assumed that $P_0=0$. Then
the monomial result gives that the above is less than or equal to

$$
\sum A(N,p)\cdot |a{_\a}|\cdot ||\n x^{\a}||_{L^p(Q_0)}
=
\sum A(N,p)\cdot ||\n a{_\a} x^{\a}||_{L^p(Q_0)}.
$$

The main trick in this proof is to use the triangle inequality for $1\le p^*$
and then be able to pass over to $p<1$ and then
use the inverse triangle inequality.
Hence the RHS above is lessthan or equal to

$$
A(N,p)\cdot (||\n P_{\oplus}||_{L^p(Q_0)}+||\n P_{\ominus}||_{L^p(Q_0)})=
A(N,p)\cdot |||\n P|||_{L^p(Q_0)}.
$$

End of proof of Theorem.
\bs

{\bf 1.4 Theorem.}
{\sl Let $1\le m$, $m$ integer. Let ${N\over N+1}\le p\le 1$.
Let $p^*={Np\over N-mp}$. For $P\in \Cal P$ and $P_{m-1}$ the 
$\Cal P_{m-1}$ part of $P$. Let $Q_0$ be as above. Then
}

$$
|||P-P_{m-1}|||_{L^{p^*}(Q_0)}
\le
A(N,p,m)\cdot|||\n^m P|||_{L^p(Q_0)}.
$$
\bs

Study the order one (say) Sobolev space with $p<1$ and $1\le p^*$.
Domain $Q_0$ as before and seminorm the previous ``triple stroke'' one.
The ordinary inequality for $p>1$ is also relied on for the inductive
proceedure.
\bs

{\bf 1.5 Definition.} 
{\sl Let $W^{1,p,|||}(Q_0)$ be the closure $\Cal P$ in norm}

$$
||u||_{W^{1,p,|||}(Q_0)}=|||u|||_{L^p(Q_0)}+|||\n u|||_{L^p(Q_0)}.
$$
\bs

Then the Sobolev imbedding for this ``triple stroke'' Sobolev space
is as follows.
\bs

{\bf 1.6 Theorem.} {\sl Let ${N\over N+1}\le p<1$ with $p^*={N\over N-p}$ 
i.e. also $1\le p^*<{N\over N-1}$ then it holds}

$$
|||u|||_{L^{p^*}(Q_0)}
\le 
A(N,p)\cdot (|||u|||_{L^p(Q_0)}+|||\n u|||_{L^p(Q_0)}).
$$
\bs

{\bf Proof.}
It is enough to study this for polynomials. Let $P\in \Cal P$.
Then by theorem 1.3 and the triangle inequality it holds

$$
|||P|||_{L^{p^*}(Q_0)}
\le
A(N,p)\cdot
(|||P_0|||_{L^{p^*}(Q_0)}
+
|||\n P|||_{L^p(Q_0)}),
$$

but

$$
|||P_0|||_{L^{p^*}(Q_0)}
=
|||P_0|||_{L^p(Q_0)}
$$
\vfill\eject

and by definition

$$
|||P_0|||_{L^p(Q_0)}
\le
|||P|||_{L^p(Q_0)}.
$$

End of proof.
\bs

In order to generalize to the order $m$ case, use a similar proceedure.
\bs

\heading
Further problems/questions
\endheading
\bs

Questions that are related either for ``triple stroke norms'' or 
in more standard way concerning polynomials possible also for
$1\le p$ as well.
\bs

I. STUDY: H\"older space in the left hand side.
\bs

II. STUDY: The question for different domains and different norm-like entities.
\bs

III. STUDY: The closure of the polynomials in Sobolev norm and any $p>0$.
Compare with closure of $C^m$-functions for a cube or general domain.
\bs

\enddocument
\bye